\documentclass{article}
\usepackage{amsmath,amssymb}

\numberwithin{equation}{section}

\newtheorem{lem}{Lemma}
\newtheorem{thm}{Theorem}
\newcommand{\itg}[0]{\mathbb{Z}}              
\newcommand{\cpx}[0]{\mathbb{C}}              
\newcommand{\la}{\lambda}
\newcommand{\ep}[0]{\varepsilon}

\newcommand{\al}[0]{\alpha}
\newcommand{\be}[0]{\beta}

\newcommand{\h}[0]{\hbar}
\newcommand{\PI}[0]{2\pi i}

\newenvironment{prf}{\noindent{\it Proof\/}.}{$\;\square$\par\medskip}

\begin{document}

\title
{Diagonalization of the elliptic Macdonald-Ruijsenaars difference system of
type $C_2$}
\author
{
Tetsuya KIKUCHI\thanks{{\it E-mail address} : tkikuchi@math.tohoku.ac.jp} 
}
\date{}

\maketitle
\centerline{{\it Mathematical Institute, Tohoku University,}}
\centerline{\it Sendai 980-8578, JAPAN}

\begin{abstract}
\noindent
We study a pair of commuting difference operators arising 
from the elliptic solution of the dynamical Yang-Baxter 
equation of type $C_2$.
The operators act on the space of meromorphic functions on 
the weight space of $\mathfrak{sp}(4,\cpx)$.
We show that these operators can be identified with the system 
by van Diejen and by Komori-Hikami with special parameters.
It turns out that our case can be related to the difference
Lam\'e operator (two-body Ruijsenaars operator)
and thereby we diagonalize the system on the finite dimensional
space spanned by the level one characters of 
the $C_2^{(1)}$-affine Lie algebra.
\end{abstract}

\section{Introduction}

The Ruijsenaars system of difference operators \cite{Ruij},
are difference analogue of the Calogero-Moser integrable 
system of differential operators.
The operators of the system is defined in terms of elliptic function, 
and in the trigonometric limit, they degenerate to the
Macdonald $q$-difference operators \cite{Macbook}.
The Ruijsenaars system has been studied extensively.
Especially, Hasegawa shows that this system can be
obtained as transfer matrices associated to the Sklyanin algebra \cite{has97}
and Felder-Varchenko reconstructed them as transfer matrices 
associated to the dynamical $R$-matrices \cite{FV2}.
These two approaches are related by the vertex-IRF 
correspondence \cite{Baxter} \cite{JMO88}.

Extending these works, in \cite{HIK} we construct
a pair of commuting difference operators 
acting on the space of functions on the $C_2$ type weight space.
The method therein is based on the elliptic solution of the dynamical 
Yang-Baxter equation of type $C_2$
(or Boltzmann weights of the $C_2^{(1)}$ face model \cite{JMO2}).
We have also shown that the space spanned by 
the level one characters of the affine Lie algebra 
$\widehat{\mathfrak{sp}}(4,\cpx)$
is invariant under the action of the difference operators.

On the other hand,
a generalization of the Ruijsenaars system 
to $BC_n$ case is studied by van Diejen \cite{vD2} and
Komori-Hikami \cite{KH97} \cite{KH98}. 
First, van Diejen constructed two elliptic commuting operators, 
one is of the first order and the other is of the $n$-th order.
Therefore he obtained an elliptic extension of difference 
Calogero-Moser system of type $BC_2$ \cite{vD1}. 
Extending this work by van Diejen,
Komori and Hikami obtained a general family
of $n$ commuting difference operators with elliptic function 
coefficients.
Besides the step parameter of difference operator 
and the modulus of elliptic functions,  
the family contains nine arbitrary parameters.
Their construction uses 
Shibukawa-Ueno's elliptic $R$-operator \cite{SU}
together with the elliptic $K$-operators \cite{KH_K-op1} \cite{KH_K-op2},
the elliptic solution to the reflection equation.
It can be regarded as an elliptic 
generalization of Dunkl type operator approach to Macdonald systems,  
which have been extensively used by Cherednik \cite{C} 
(see \cite{N} for $BC_n$ case). 

This paper has two goals.
One is to establish the relationship between our system of difference operators 
and the van Diejen-Komori-Hikami system.
The other is to diagonalize our difference operators on 
the finite dimensional space spanned by theta functions.
The first goal is attained in section \ref{taiou} and
second is in section \ref{diag}.

In section \ref{taiou} we review the construction of the
elliptic difference system of type $C_2$ and give a new form
of our operators.
After this, we will establish an identity consisting 
of theta functions (Lemma \ref{new}),
and explain how our system can be identified with
from van Diejen-Komori-Hikami system
with special choice of parameters (Theorem \ref{connect}).
That is, our approach to the difference operators as
transfer matrices, based on the knowledge of the Boltzmann
weights, reproduces a special case among the family of commuting operators
obtained by Dunkl-type approach.
It should be also mentioned that those two approaches to the system
are not yet related, although the resulting commuting
operators have the relationship as above.

In section \ref{diag}, we introduce the finite dimensional space
of theta functions invariant under the action of 
Weyl group and its basis after Kac-Peterson \cite{KP}. 
Our aim is to diagonalize our operators on this space (Theorem \ref{main}).
This is an elliptic analogue of the eigenvalue problem of 
Macdonald operators 
on the space of symmetric polynomials.
Their eigenfunctions, called Macdonald-Koornwinder polynomials,
are much investigated  in $q$-orthogonal 
polynomial theory \cite{Koorn} \cite{N}.

\section{The difference operators of type $C_2$}\label{taiou}

\subsection{Construction of the difference operators  of type $C_2$}

Let $\mathfrak{g}$ be the Lie algebra $\mathfrak{sp}(4,\cpx)$, 
$\mathfrak{h}$ its Cartan subalgebra and
$\mathfrak{h}^*$ the dual space of $\mathfrak{h}$. 
We realize the root system $R$ for $(\mathfrak{g},\mathfrak{h})$
as $R:=\{ \pm(\ep_1\pm\ep_2), \pm 2\ep_1, \pm 2\ep_2\}\subset \mathfrak{h}^*.$
A normalized Killing form $(\,,\,)$ is given by
\begin{equation}\label{killing}
(\ep_j,\ep_k)=\frac{1}{2}\delta_{jk},
\end{equation}
and the square length of the long roots $\pm 2 \ep_i$ is two.
We will identify the space $\mathfrak{h}$ and its
dual $\mathfrak{h}^*$ via the form $(\,,\,)$.
The fundamental weights
are given by $\Lambda_1=\ep_1,\Lambda_2=\ep_1+\ep_2$.
Let ${\cal{P}}_d$ be the set of weights for the fundamental representation
$L(\Lambda_d)$.
We have
\begin{equation}
{\cal{P}}_1=\{\pm\ep_1,\pm\ep_2\},\quad
{\cal{P}}_2=\{\pm(\ep_1 \pm \ep_2), 0 \}.\label{Pd}
\end{equation}

Let $d,d'$ be $1$ or $2.$
The $C_2^{(1)}$ type Boltzmann weights of type $(d,d')$
are given as follows. Fix a complex parameter $\h \in \cpx$.
For any pair $\la, \mu, \nu, \kappa \in \mathfrak{h}^*$ 
of weights, the Boltzmann weight
\[ 
W_{dd'}\!
\left(\left.\begin{array}{ll}
        \la     &       \mu     \\
        \kappa  &       \nu     \\
\end{array}\,\right|u
\right)
\]
is a function of the spectral 
parameter $u \in \cpx$ and $\la \in \mathfrak{h}$.
They satisfy the condition
\begin{equation*}
W_{dd'}\!
\left(\left.\begin{array}{ll}
        \la     &       \mu     \\
        \kappa  &       \nu     \\
\end{array}\,\right|u
\right)=0
\;\mbox{unless}\;\mu-\la,\nu-\kappa\in 2\h{\cal P}_d,\;
\kappa-\la,\nu-\mu \in 2\h{\cal P}_{d'},
\end{equation*}
and solve the Yang-Baxter equation of the face type,
\begin{align}
&\sum_\eta W_{dd'}
\left(\left.\begin{array}{ll}
        \rho    &       \eta    \\
        \sigma  &       \kappa  \\
\end{array}\,\right|u-v
\right)
W_{dd''}
\left(\left.\begin{array}{ll}
        \la     &       \mu     \\
        \rho    &       \eta    \\
\end{array}\,\right|u-w
\right)
W_{d'd''}
\left(\left.\begin{array}{ll}
        \mu     &       \nu     \\
        \eta    &       \kappa  \\
\end{array}\,\right|v-w
\right)\nonumber\\
=&
\sum_\eta W_{d'd''}
\left(\left.\begin{array}{ll}
        \la     &       \eta    \\
        \rho    &       \sigma  \\
\end{array}\,\right|v-w
\right)
W_{dd''}
\left(\left.\begin{array}{ll}
        \eta    &       \nu     \\
        \sigma  &       \kappa  \\
\end{array}\,\right|u-w
\right)
W_{dd'}
\left(\left.\begin{array}{ll}
        \la     &       \mu     \\
        \eta    &       \nu     \\
\end{array}\,\right|u-v
\right).
\label{YBE1}
\end{align}
This equation is also known as the dynamical Yang-Baxter equation.
Here we give the explicit formula for $W_{11}$
and see \cite{HIK} for the other type $W_{dd'}$ ($(d,d') = (1,2), \, (2,1), \, (2,2)$)
which are obtained by fusion procedure.
They are expressed by the Jacobi theta function
$\theta_1(u) = \theta_1(u|\tau)$ 
with elliptic modulus $\tau$ in the upper half plane $\mathfrak{H}_+$
(See Appendix B for the definition of $\theta_1(u)$).
For $p,q,r,s\in{\cal P}$ such that $p+q=r+s$, we will write 
\begin{equation*}
\begin{matrix}    
                &p              &               \\
   s\!\!\!\!    &\boxed{u}      &\!\!\!\!\! q   \\
                &r              &               \\
\end{matrix}
=
W_{11} \left( \left.\begin{array}{ll}
          \la   &        \la+ 2\h p        \\
         \la+ 2\h s &\la+ 2\h (p + q)    \\
\end{array}
\;\right|u
\right).
\label{face}
\end{equation*}
The explicit formula for $W_{11}$ is given as follows:
\begin{align}
\begin{matrix}    
                &p       &      \\
   p\!\!\!\!    &\boxed{u} &\!\!\!\!\! p\\
                &p      &       \\
\end{matrix}&=
\frac{\theta_1(c-u)\,\theta_1(u+\h)}
{\theta_1(c)\,\theta_1(\h)},\label{2-5a}
\\
\begin{matrix}  
                &p       &      \\
          p \!\!\!\!    &\boxed{u} &\!\!\!\!\! q\\
                &q      &       \\
\end{matrix}&=
\frac{\theta_1(c-u)\,\theta_1(\la_{p-q}-u)}
{\theta_1(c)\,\theta_1(\la_{p-q})} \qquad (p\ne \pm q),\label{2-5b}
\\
\begin{matrix}    
                &q       &      \\
          p \!\!\!\!     &\boxed{u} &\!\!\!\!\! p\\
                &q      &       \\
\end{matrix}&=
\frac{\theta_1(c-u)\,\theta_1(u)\,\theta_1(\la_{p-q}+\h)}
{\theta_1(c)\,\theta_1(\h)\,\theta_1(\la_{p-q})} 
\qquad (p\ne \pm q),\label{2-5c}
\\
\begin{matrix}    
                &q       &      \\
          p \!\!\!\!    &\boxed{u} &\!\!\!\!\! -q\\
                &-p &       \\
\end{matrix}&=
-\frac{\theta_1(u)\,\theta_1(\la_{p+q}+\h+c-u)}
{\theta_1(c)\,\theta_1(\la_{p+q}+\h)} 
\frac{\theta_1(2\la_p+2\h)}{\theta_1(2\la_q)}
\frac{\prod_{r\ne \pm p}\theta_1(\la_{p+r}+\h)}
{\prod_{r\ne \pm q}\theta_1(\la_{q+r})} \qquad (p\ne q),\label{2-5d}
\\
\begin{matrix}    
                &p       &      \\
          p \!\!\!\!    &\boxed{u} &\!\!\!\!\! -p\\
                &-p &       \\
\end{matrix}&=
\frac{\theta_1(c-u)\,\theta_1(2\la_p+\h-u)}
{\theta_1(c)\,\theta_1(2\la_p+\h)} \\ \nonumber
&-\frac{\theta_1(u)\,\theta_1(2\la_p+\h+c-u)}
{\theta_1(c)\,\theta_1(2\la_p+\h)}
\frac{\theta_1(2\la_p+2\h)}{\theta_1(2\la_p)}
\prod_{q\ne\pm p}
\frac{\theta_1(\la_{p+q}+\h)}{\theta_1(\la_{p+q})}.\label{2-5e}
\end{align}
Here the crossing parameter $c$ is fixed to be $c:=-3\h$.

We define the difference operators $M_d(u)\;(u\in \cpx,d=1,2)$ acting on 
the functions on $\mathfrak{h}$ by means of the Boltzmann weights of
type $(1,2)$ and $(2,2)$.
\begin{equation*}
(M_d(u)f)(\la):=\sum_{p\in {\cal{P}}_d}
{W_{d2}}\left(\left.\begin{array}{ll}    
          \la   &        \la+2\h{p}        \\
          \la   &        \la+2\h{p}       
\end{array}\right|u \right)
\;T^{\h}_{2p} f(\la).
\end{equation*}
Here the shift operator $T_{2p}^\h$ is defined as
$$
T^\h_{2p} f(\la) : =f(\la+2\h{p}\,).
$$
For $\la\in\mathfrak{h}^*$ and $p\in{\cal P}_d\,(d=1,2)$,
we put
\begin{equation*}
\la_{p}:=(\la,p)
  \label{lambda_p}.
\end{equation*} 
Note that if we denote $\la_i = (\la,\ep_i) \; (i=1,2)$ and $f(\la) = f(\la_1,\la_2)$,
then 
$$T^\h_{\pm 2\ep_1} f(\la_1,\la_2) =f(\la_1 \pm \h,\la_2), \quad
T^\h_{\pm 2\ep_2} f(\la_1,\la_2) =f(\la_1,\la_2 \pm \h).
$$

\begin{thm}\label{commute}
{\rm \cite{HIK}}

$({\rm i})$
For each $u,v \in \cpx$, we have 
$M_d(u)M_{d'}(v) = M_{d'}(v)M_d(u) \; (d,d'=1,2).$

$({\rm ii})$
The explicit form of $M_d(u)$ are as follows $:$
\begin{equation}\label{M1}
M_1(u) = F(u) \,
\sum_{p\in{\cal{P}}_1}\prod_{\substack{q\in{\cal{P}}_1 \\ q\ne \pm p}}
\frac{\theta_1(\la_{p+q}-\h)}{\theta_1(\la_{p+q})}
T^\h_{2p},
\end{equation}
\begin{equation}\label{M2}
M_2(u) = G(u) \left(
\sum_{\substack{p=\pm\ep_1 \\ q=\pm\ep_2}}
\left(\frac{\theta_1(\la_{p+q}-\h)}{\theta_1(\la_{p+q}+\h)}
T^\h_{2p} T^\h_{2q}
+ U(\la_p,\la_q)
\right) -H(u) \right).
\end{equation}
Here $U(\la_p,\la_q)$ is given by $:$
\[
U(\la_p,\la_q) = 
\frac{\theta_1(2\h)}{\theta_1(6\h)}
\frac{\theta_1(2\la_p+2\h)\,\theta_1(2\la_q+2\h)}
{\theta_1(2\la_p)\,\theta_1(2\la_q)}
\frac{\theta_1(\la_{p+q}-5\h)\,\theta_1(\la_{p+q}+2\h)}
{\theta_1(\la_{p+q})\,\theta_1(\la_{p+q}+\h)},
\]
and $F(u),G(u),H(u)$ are the following functions depend only on $u$ and $\h :$
\[ F(u) :=  \frac{\theta_1(u)\,\theta_1(u+2\h)^2\,\theta_1(u+4\h)}
{\theta_1(-3\h)^2\,\theta_1(\h)^2} , \]
\begin{equation*}
 G(u) := \frac{\theta_1(u-\h)\,\theta_1(u)^2\,\theta_1(u+\h)\,
 \theta_1(u+2\h)\,\theta_1(u+3\h)^2\,\theta_1(u+4\h)}
 {\theta_1(-3\h)^4\,\theta_1(\h)^4},
\end{equation*}
\[ \mbox{and}\qquad H(u) :=
\frac{\theta_1(u+6\h)\,\theta_1(u-3\h)\,\theta_1(2\h)}
{\theta_1(u)\,\theta_1(u+3\h)\,\theta_1(6\h)}. \]
\end{thm}

The following Lemma is the key for the identification
with van Diejen's system as well as for the diagonalization of 
our difference operators. The author is grateful to van Diejen
for the information.
\begin{lem}\label{new}
We have
\begin{equation}\label{Diejen}
\sum_{\substack{p=\pm \ep_1 \\ q=\pm \ep_2}} 
U(\la_p,\la_q) 
- \sum_{\substack{p=\pm \ep_1 \\ q=\pm \ep_2}} 
\frac{\theta_1(\la_{p+q}-\h)\,\theta_1(\la_{p+q}+2\h)}
{\theta_1(\la_{p+q})\, \theta_1(\la_{p+q}+\h)} = K,
\end{equation}
where K is a constant given by
\begin{equation*}
K = \frac{\theta_1(8\h)\,\theta_1(\h)}{\theta_1(6\h)\,\theta_1(5\h)}
+ \frac{\theta_1(5\h)\,\theta_1(2\h)}{\theta_1(4\h)\,\theta_1(3\h)}
+ \frac{\theta_1(6\h)\,\theta_1(3\h)}{\theta_1(5\h)\,\theta_1(4\h)}
+ \frac{\theta_1(4\h)\,\theta_1(\h)}{\theta_1(3\h)\,\theta_1(2\h)}.
\end{equation*}
\end{lem}
\begin{prf}
Let $f(\la_p)$ be the left-hand side
of (\ref{Diejen}), regarded as a function of $\la_p \; (p \in I)$.
It is doubly periodic function of the periods $1,\tau$.
Let us show that it is entire.
The apparent poles of $f(\la_p)$ are located at
$$ \la_p=\la_q,\; \la_p=\la_q - \h \;
(p,q\in I,\; p+q \ne 0),\;\la_p=0\,(p \in I).$$
Note that $f(\la_p)$ is $W$-invariant, then 
the points $\la_p=\la_q$ and $\la_p=0$ are regular.
Also, the residue of $f(\la_p)$ at $\la_p= -\la_q- \h$ is
\[ \frac{\theta_1(2\h)}{\theta_1(6\h)} 
\frac{\theta_1(-2\la_q)\,\theta_1(2\la_q+2\h)}
{\theta_1(-2\la_q-2\h)\,\theta_1(2\la_q)}
\frac{\theta_1(-6\h)\,\theta_1(\h)}{\theta_1(-\h)} 
- \frac{\theta_1(-2h)\,\theta_1(\h)}{\theta_1(-\h)} = 0.\]

Now we have proved that $f(\la_p)$ is independent of $\la_p$, then 
we consider $g(\la_q)=f(-\la_q-2\h)$ as a function of $\la_q \;(q\ne p \in I)$:
\begin{eqnarray*}
g(\la_q) &=&\frac{\theta_1(2\h)}{\theta_1(6\h)} \left( 
\frac{\theta_1(2\la_q+2\h)\,\theta_1(2\la_q-2\h)\,\theta_1(2\la_q+7\h)}
{\theta_1(2\la_q+4\h)\,\theta_1(2\la_q+2\h)\,\theta_1(2\la_q+\h)} \right. \\
&&+ \frac{\theta_1(2\la_q+6\h)\,\theta_1(2\la_q+2\h)\,\theta_1(2\la_q-3\h)}
{\theta_1(2\la_q)\,\theta_1(2\la_q+2\h)\,\theta_1(2\la_q+3\h)}  
\\
&&+ \left. \frac{\theta_1(2\la_q+6\h)\,\theta_1(2\la_q-2\h)\,\theta_1(-3\h)\,\theta_1(4\h)}
{\theta_1(2\la_q+4\h)\,\theta_1(2\la_q)\,\theta_1(2\h)\,\theta_1(3\h)}
\right) \\
&& -\frac{\theta_1(-2\la_q-3\h)\,\theta_1(-2\la_q)}
{\theta_1(-2\la_q-2\h)\,\theta_1(-2\la_q-\h)} 
- \frac{\theta_1(2\la_q+\h)\,\theta_1(2\la_q+4\h)}
{\theta_1(2\la_q+2\h)\,\theta_1(2\la_q+3\h)}
-\frac{\theta_1(\h)\,\theta_1(4\h)}{\theta_1(2\h)\,\theta_1(3\h)}. 
\end{eqnarray*}
By the same argument we can show that $g(\la_q)$ is independent of $\la_q$.
Therefore we get $K$ by putting $\la_q=\h$ in $g(\la_q)$ 
and the proof completes.
\end{prf}

\subsection{Identification with van Diejen's system}

We define the difference operators $\widetilde{M_d}$ to be 
the components of $M_d(u)$ independent of $u$:
\begin{equation}
\widetilde{M_1} = 
\sum_{p\in{\cal{P}}_1}\prod_{\substack{q\in{\cal{P}}_1 \\ q\ne \pm p}}
\frac{\theta_1(\la_{p+q}-\h)}{\theta_1(\la_{p+q})}
T^\h_{2p},
\end{equation}
\begin{equation}
\widetilde{M_2} =  
\sum_{\substack{p=\pm\ep_1 \\ q=\pm\ep_2}}
\left(\frac{\theta_1(\la_{p+q}-\h)}{\theta_1(\la_{p+q}+\h)}
T^\h_{2p} T^\h_{2q}
+
\frac{\theta_1(\la_{p+q}-\h)\,\theta_1(\la_{p+q}+2\h)}
{\theta_1(\la_{p+q})\,\theta_1(\la_{p+q}+\h)}
\right) 
\end{equation}

More general commuting
difference operators ${\cal H}_1, {\cal H}_2$ are obtained by 
van Diejen and later by Komori-Hikami in a different way.
In this subsection we identify our operators $\widetilde{M_1}, \widetilde{M_2}$ as
van Diejen's system of difference operators with special values of parameters.
The operators ${\cal H}_1, {\cal H}_2$ depend on nine 
complex parameters $\mu, \mu_r,\mu_r' \; (r=0,1,2,3)$
satisfying the condition 
\begin{equation}\label{zero-wa}
\sum_r  \left( \mu_r + \mu_r' \right)=0
\end{equation}
and are defined by
\begin{eqnarray*}
{\cal H}_1 &=&  
\sum_{\ep = \pm 1}
w(\ep x_1)v(\ep x_1+ x_2)v(\ep x_1- x_2) T_{\ep 1}^\gamma  \\
&& {} + \sum_{\ep = \pm 1}
w(\ep x_2)v(\ep x_2+x_1)v(\ep x_2-x_1) T_{\ep 2}^\gamma  \; + \;
U_{\{1,2\},1}, 
\end{eqnarray*}
\begin{eqnarray*}
{\cal H}_2 &=& 
\sum_{\ep,\ep ' = \pm 1}
w(\ep x_1) w(\ep ' x_2) 
v(\ep x_1+ \ep ' x_2)v(\ep x_1 + \ep 'x_2 + \gamma) 
T_{\ep 1}^\gamma  T_{\ep '2}^\gamma \\
&& {} + U_{\{2 \},1} \sum_{\ep = \pm 1} 
w(\ep x_1)v(\ep x_1+x_2)v(\ep x_1-x_2) T_{\ep 1}^\gamma  \\
&& {} + U_{\{1 \},1} \sum_{\ep = \pm 1} 
w(\ep x_2)v(\ep x_2+x_1)v(\ep x_2-x_1) T_{\ep 2}^\gamma \; + \; 
U_{\{1,2\},2}. 
\end{eqnarray*}
Here $T^\gamma_{\pm i} \; (i=1,2)$ stand for the shift operators 
\[ T^\gamma_{\pm1}f(x_1,x_2) = f(x_1 \pm \gamma,x_2), \quad
T^\gamma_{\pm2}f(x_1,x_2) = f(x_1,x_2 \pm \gamma) \]
and
\begin{equation}\label{vz}
v(z) := \frac{\sigma(z+\mu)}{\sigma(z)}, \quad
w(z) := \prod_{0 \le r \le 3}
\frac{\sigma_r(z+\mu_r) \, \sigma_r(z+\mu'_r+\gamma/2)}
{\sigma_r(z)\,\sigma_r(z+\gamma/2)}, 
\end{equation}
where $\sigma(z) = \sigma_0(z)$ denotes the sigma function with two quasi 
periods $\omega_1, \; \omega_2$ and $\sigma_r(z) \; (r=1,2,3)$ associated 
function obtained by shift of argument over the half periods 
(See Appendix B  for more detail). The functions $U_{\{j\},1},  
U_{\{ 1,2 \},j} \; (j=1,2)$ are defined as follows :
\[ U_{\{j\},1} = -w(x_j) -w(-x_j) \quad (j=1,2), \]
\[  U_{\{ 1,2 \},1}
= \sum_{0 \le r \le 3} c_r \prod_{j=1,2}
\frac{\sigma_r(\mu - \gamma/2 +x_j)\,\sigma_r(\mu - \gamma/2 - x_j)}
{\sigma_r(-\gamma/2 +x_j)\,\sigma_r(- \gamma/2 - x_j)},
\]
where
\[ c_r = \frac{2}{\sigma(\mu)\,\sigma(\mu-\gamma)}
\prod_{0 \le s \le 3}\sigma_s(\mu_{\pi_r(s)}-\gamma/2)\,
\sigma_s(\mu'_{\pi_r(s)}), \]
with $\pi_r$ denoting the permutation $\pi_0 = id$, $\pi_1 = (01)(23)$,
$\pi_2 = (02)(13)$, $\pi_3 = (03)(12)$.
\begin{equation}\label{U122}
U_{\{ 1,2 \},2}
= \sum_{\ep,\ep ' \in \{ 1,-1 \}}
w(\ep x_1) w(\ep' x_2) 
v(\ep x_1+ \ep' x_2) v(-\ep x_1 - \ep'x_2 - \gamma) 
\end{equation}

We mention that the Komori-Hikami system in \cite{KH98} is of more complicated
form and has nine arbitrary parameters, that is, they removed 
the condition (\ref{zero-wa}).

In ${\cal H}_1,{\cal H}_2$, we specialize parameters $\mu,\mu_r,\mu'_r \;
(r = 0,1,2,3)$ as $\mu = -\gamma,\,\mu_r = \mu'_r = 0$. 
Then $w(z) = 1$ and $U_{\{ 1,2 \} , 1} = 0$.
Let us denote these specialized operators by $\bar{\cal H}_1,
\bar{\cal H}_2$.
Because of these simplifications, we immediately obtain the 
following from Lemma \ref{new}, giving the identification of our
system $\{ \widetilde{M_1}, \; \widetilde{M_2} \}$ and van 
Diejen's $\{ \bar{\cal H}_1, \; \bar{\cal H}_2 \}$.

\begin{thm}\label{connect} 
For a function $f(\la_1,\la_2)$
on  $\mathfrak{h}$, we set $\varphi(f)(x_1,x_2)$ by
$$ \varphi(f)(x_1,x_2) :=  \exp  \frac{\eta_1(x_1^2 + x_2^2)}{\omega_1} 
f(\frac{x_1}{2\omega_1},\frac{x_2}{2\omega_1}), $$
and  let $\gamma = 2\omega_1 \h$, we have
\begin{align*}
\varphi \; \widetilde{M_1} \; \varphi^{-1}  &= 
e^{2\eta_1 \gamma^2/\omega_1} \,
\bar{\cal H}_1 , \\
\varphi \; \widetilde{M_2} \; \varphi^{-1} &= 
e^{2\eta_1 \gamma^2/\omega_1} \,
\left( \bar{\cal H}_2 + 2 \bar{\cal H}_1 \right). 
\end{align*}
\end{thm}
\begin{prf}
Use the connection between the theta function and sigma function (\ref{tesig}) in 
Appendix B and (\ref{Diejen}) to 
compare (\ref{M2}) and (\ref{U122}).
\end{prf}


\section{Diagonalization of the system}\label{diag}

\subsection{The space of theta functions}

Let $Q$ and $Q^\vee$ be the root and coroot 
lattice, $P$ and $P^\vee$  the weight and coweight lattice respectively.
Under the identification $\mathfrak{h}=\mathfrak{h}^*$ via the form $(\,,\,),$
they are given by 
\begin{equation}\label{root}
P= \sum_{j=1,2}\itg \ep_j, 
\quad
Q^\vee=\sum_{j=1,2}\itg 2\ep_j,
\end{equation}
and
\[
P^\vee=  \; Q = \itg 2\ep_1  + \itg 2\ep_2  + \itg (\ep_1 + \ep_2).
\]
For $\beta\in {\mathfrak h}^*$, we introduce the following 
operators $T_{\tau\beta},T_{\beta}$ acting on the functions on
$\mathfrak{h}^*$: 
\begin{align*}
(T_{\be}f)(\la) &:=f(\la+\be), \\
(T_{\tau\be}f)(\la)
&:=\exp\left[\PI \left((\la,\be)+\frac{(\be,\be)}{2}\tau \right)\right]
f(\la+\tau\be)
\end{align*}
We define the space of theta functions (of level $1$) by
\begin{equation*}
Th_1
:=\left\{ f \;\mbox{is holomorphic on} \; \mathfrak{h}^* \; |
\; T_{\tau\al}f=T_\al f=f  \quad (\forall\al\in Q^\vee) 
 \right\}.
\end{equation*}
For each $\mu \in P $ and fixed 
$\tau \in \mathfrak{H}_+$, 
we define the classical theta function $\Theta_{\mu}(\la)$ of 
$\la \in  \mathfrak{h}^* $ by
\[ \Theta_{\mu}(\la) :=
\sum_{\gamma \in \mu + Q^\vee}  \exp \left[ 2\pi i 
\left( (\gamma, \la) + \frac{(\gamma,\gamma)}{2} \tau \right)  \right] .
\]
It is known that
\[ \{ \Theta_{\mu}(\la) \; | \; \mu \equiv 0, \, \ep_1, \, \ep_2, \, 
\ep_1+\ep_2 \; \mathrm{mod}\, Q^\vee \} \]
gives a basis for $Th_1$ over $\cpx$ \cite{KP}.

Let $W \subset GL(\mathfrak{h}^*) $ denote the Weyl group for
$(\mathfrak{g},\mathfrak{h})$, and consider the $W$-invariants in $Th_1$: 
\begin{equation*}
Th^W_1 := \left\{ f \in Th_1 \; | \;
f(w\la)=f(\la) \; (\forall w\in W)
\right\}.
\end{equation*}
\begin{thm}{\rm \cite{HIK}}
The operators $\widetilde{M_1}$, $\widetilde{M_2}$
preserves $Th_1^W$.
\end{thm}
For $\mu \in P$ , we define $W_\mu := \{ w \in W \; | \; w \mu = \mu \}$ 
and introduce the following symmetric sum of 
theta functions, 
\[ S_{\mu}(\la) := \frac{1}{|W_\mu|}
\sum_{w \in W} \Theta_{w(\mu)}(\la). \]
Then
$$
\{ S_\mu(\la) \; | \; \mu \equiv \, 0, \, \Lambda_1 \,(=\ep_1), \, 
\Lambda_2 \,(=\ep_1+\ep_2) \;
\mathrm{mod} \, Q^\vee \} 
$$ 
forms a basis for $Th_1^W$ over $\cpx$.

It is known that $Th_1^W$ is also spanned by the 
level $1$ characters of the affine Lie 
algebra $\widehat{\mathfrak{sp}}(4,\cpx)$. 
Note that $\Theta_{-\mu}(\la) = \Theta_\mu(\la)$ and
$\Theta_{\ep_1+\ep_2}(\la) = \Theta_{\ep_1-\ep_2}(\la)$.
So that we have
$$ 
S_0(\la) = \Theta_0(\la), \quad S_{\Lambda_1}(\la) 
= 2( \Theta_{\ep_1}(\la) + \Theta_{\ep_2}(\la) ), \quad
S_{\Lambda_2}(\la) = 4\Theta_{\ep_1+\ep_2}(\la).
$$

\subsection{Diagonalization of $\widetilde{M_d}$ }

In this subsection, we diagonalize  the
operators $\widetilde{M_d}$ on the space $Th_1^W$. 
We set 
\begin{align*}
f_1(\la) &:= \Theta_{\ep_1}(\la) + \Theta_{\ep_2}(\la), \quad
f_2(\la) := \Theta_0(\la) + \Theta_{\ep_1+\ep_2}(\la) \quad \mbox{and} \\
f_3(\la) &:= \Theta_0(\la) - \Theta_{\ep_1+\ep_2}(\la). 
\end{align*}
They are linearly independent in the space $Th^W_1$.

\begin{thm}\label{main}
The functions $f_i(\la)$ $(i = 1,2,3)$ are 
common eigenfunctions of $\widetilde{M_d}:$
$$ \widetilde{M_d} f_i(\la) = E_{d,i} f_i(\la) \quad (d=1,2, \; i = 1,2,3).$$
The eigenvalues are given by
$$ E_{1,i} = 
\left( \frac{\theta_1(2\h) \theta_{i+1}(0) }
{\theta_1(\h) \theta_{i+1}(\h)} \right)^2 $$
and $E_{2,i} = 2 E_{1,i}$,
where the Jacobi theta 
functions $\theta_i(z) = \theta_i(z|\tau) \; (i=2,3,4)$  are 
defined as in Appendix B.
\end{thm}

We will prove this theorem by using the following three lemmas.
First,  we show that the operators $\widetilde{M_d}$ split 
into two $A_1$-type components.
\begin{lem}\label{split}
Let us denote $\la_\pm := (\la,\ep_1\pm \ep_2)$ 
and define
\[ H_\pm := 
\frac{\theta_1(\la_\pm -\h)}{\theta_1(\la_\pm)}
T^\h_{\ep_1\pm \ep_2}
+\frac{\theta_1(-\la_\pm -\h)}{\theta_1(-\la_\pm)}
T^\h_{-(\ep_1\pm \ep_2)}.
\]
Then we have
\begin{equation}\label{takemura}
\widetilde{M_1} = H_+H_-, \quad
\widetilde{M_2} = H_+^2 + H_-^2.  
\end{equation}
\end{lem}
\begin{prf}
To prove the first identity, we note that
\begin{align*}
& \frac{\theta_1(\la_+ -\h)}{\theta_1(\la_+ )}
T^\h_{\ep_1+\ep_2}
\frac{\theta_1(\la_- -\h)}{\theta_1(\la_-)}
T^\h_{\ep_1-\ep_2} \\
&= \frac{\theta_1(\la_+ -\h)}{\theta_1(\la_+)}
\frac{\theta_1((\la+\h(\ep_1+\ep_2))_- -\h)}
{\theta_1((\la+\h(\ep_1+\ep_2))_-)}
T^\h_{\ep_1+\ep_2} T^\h_{\ep_1-\ep_2} \\
&= \frac{\theta_1(\la_+ -\h)}{\theta_1(\la_+ )}
\frac{\theta_1(\la_- -\h)}
{\theta_1(\la_-)}
T^\h_{2\ep_1}.
\end{align*}
Here we used the identity $(\ep_1 + \ep_2, \ep_1 - \ep_2) = 0$.
The second identity follows from, for instance,  
\begin{align*}
& \frac{\theta_1(\la_+ -\h)}{\theta_1(\la_+)}
T^\h_{\ep_1+\ep_2}
\frac{\theta_1(\la_+ -\h)}{\theta_1(\la_+)}
T^\h_{\ep_1+\ep_2} \\
&= \frac{\theta_1(\la_+ -\h)}{\theta_1(\la_+)}
\frac{\theta_1((\la+\h(\ep_1+\ep_2))_+ -\h)}
{\theta_1((\la+\h(\ep_1+\ep_2))_+)}
T^\h_{\ep_1+\ep_2} T^\h_{\ep_1+\ep_2} \\
&= \frac{\theta_1(\la_+ -\h)}{\theta_1(\la_+)}
\frac{\theta_1(\la_+ + \h -\h)}
{\theta_1(\la_+ + \h)}
T^\h_{2(\ep_1+\ep_2)} \\
&= \frac{\theta_1(\la_+ -\h)}
{\theta_1(\la_+ + \h)}
T^\h_{2\ep_1} T^\h_{2\ep_2}.
\end{align*}
Here we used the identity $(\ep_1 + \ep_2, \ep_1 + \ep_2) = 1$.
\end{prf}

Second, we consider the eigenvalue problem 
for the $A_1$-type difference operator 
(difference Lam\'e or two-body Ruijsenaars operator)
\begin{equation}\label{mondai}
\frac{\theta_1(z-\ell\h)}{\theta_1(z)} f(z+\h) + 
\frac{\theta_1(z+\ell\h)}{\theta_1(z)} f(z-\h) = 
E f(z) .
\end{equation}
 
\begin{lem} 
For the special coupling constant $\ell = 1$, the functions
\[ \theta_i(z) \quad (i = 2,3,4) \]
are solutions of the equation  \rm{(\ref{mondai})} with eigenvalues
\[ E = E_i = 
\frac{\theta_1(2\h)\,\theta_i(0)}{\theta_1(\h)\,\theta_i(\h)} 
\quad (i=2,3,4). \]
\label{a1eigen}
\end{lem} 
\begin{prf}
We note that the functions $\theta_2(z), \theta_3(z)$, and $\theta_4(z)$ can
be rewritten in $e^{cz} \theta_1(z+t)$ up to a constant, where
\begin{equation}
(t,c) = (\frac{1}{2}, 0), \; (\frac{1+\tau}{2}, \pi i) \; \; \mbox{and} \; \;
(\frac{\tau}{2}, \pi i)
\label{BAsol}
\end{equation}
respectively (See the formula (\ref{theta2}) in Appendix B).
Let $(t,c)$ be one of these, and
we denote by $g(z)$ the function which obtained by the 
action of the $A_1$-type difference 
operator (\ref{mondai}) with $\ell = 1$ to $e^{cz} \theta_1(z+t)$:  
\begin{equation*}
g(z) := \frac{\theta_1(z-\h)\,\theta_1(z+t+\h)}{\theta_1(z)} e^{c(z+\h)} +
\frac{\theta_1(z+\h)\,\theta_1(z+t-\h)}{\theta_1(z)} e^{c(z-\h)}.
\end{equation*}
This is holomorphic and doubly quasi-periodic function: 
$$g(z+1) 
= -e^c g(z), \quad g(z+\tau) = -e^{\pi i \tau - 2\pi i (z+t) + c\tau} g(z).
$$ 
Moreover, $g(z)=0$ at $z=-t$.
Therefore, $g(z)$ is equal to $e^{cz} \theta_1(z+t)$ up to a constant, 
which is the value of 
\begin{equation}
\frac{\theta_1(z-\h)\,\theta_1(z+t+\h)}{\theta_1(z)\, \theta_1(z+t)} e^{c\h} +
\frac{\theta_1(z+\h)\,\theta_1(z+t-\h)}{\theta_1(z)\, \theta_1(z+t)} e^{-c\h}
\label{bethe}
\end{equation}
at any chosen point.
If we choose $z=\h$, then the first term in (\ref{bethe}) vanishes and we have
$$
\frac{\theta_1(2\h)\,\theta_1(t)}{\theta_1(\h)\, \theta_1(\h+t)} e^{-c\h}
= \frac{\theta_1(2\h)\,\theta_i(0)}{\theta_1(\h)\, \theta_i(\h)},
$$
where $i =2,3$ and $4$ corresponding to the 
values of $(t,c)$ in (\ref{BAsol}), as an eigenvalue.
\end{prf}

{\bf Remark}.
This can be regarded as a special case of Felder-Varchenko's study \cite{FV1}. 
They expressed the solutions of (\ref{mondai}) in terms of 
the algebraic Bethe Ansatz method, which is originally developed and
applied to the spin chain model. In fact, the operator in 
the left hand side of (\ref{mondai}) can be regarded 
as the transfer matrix of the simplest spin chain, that is, it
consists of only one site of freedom with spin $\ell = 1$.
In this case, the Bethe Ansatz equation 
\begin{equation}
\frac{\theta_1(t-\h)}{\theta_1(t+\h)} = e^{2\h c}
\end{equation}
is exactly the same as the condition that the function (\ref{bethe}) dose not 
have pole at $z=-t$.

Because of Lemma \ref{a1eigen}, the product of the theta functions 
\[ \theta_i(\la_+) \theta_j(\la_-)
\quad (i,j = 2,3,4)
\]
are simultaneous eigenfunctions of the 
operators $H_+^2,H_-^2$ and $H_+H_-$ with eigenvalues
\[
\frac{\theta_1(2\h)^2 \, \theta_i(0)^2}
{\theta_1(\h)^2 \, \theta_i(\h)^2}, \quad 
\frac{\theta_1(2\h)^2 \, \theta_j(0)^2}
{\theta_1(\h)^2 \, \theta_j(\h)^2} \quad
\mbox{and} \quad 
\frac{\theta_1(2\h)^2 \, \theta_i(0) \, \theta_j(0)}
{\theta_1(\h)^2 \, \theta_i(\h) \, \theta_j(\h)} 
\]
respectively.
Finally, we shall establish the
relationship of these Bethe Ansatz solutions and
the bases of $Th_1^W$.
\begin{lem}\label{kitei}
The functions $f_i(\la) \in Th^W_1$ are expressed 
in terms of the Jacobi theta functions as follows$:$
\begin{equation*}
f_1(\la) = 
\theta_2(\la_+) \theta_2(\la_-), \quad
f_2(\la) = 
\theta_3(\la_+) \theta_3(\la_-), \quad
f_3(\la) = 
\theta_4(\la_+) \theta_4(\la_-).
\end{equation*}
\end{lem}
\begin{prf}
Because of  the definitions of coroot lattice $Q^\vee$ (\ref{root}) 
and Killing form (\ref{killing}),
each basis of $Th_1$ is expressed as
\begin{align*}
\Theta_0(\la) &= \theta_3(2\la_1|2\tau)\theta_3(2\la_2|2\tau),  \\
\Theta_{\ep_1}(\la) &= \theta_3(2\la_1|2\tau)\theta_2(2\la_2|2\tau),  \\
\Theta_{\ep_2}(\la) &= \theta_2(2\la_1|2\tau)\theta_3(2\la_2|2\tau),  \\
\Theta_{\ep_1+\ep_2}(\la) &= \theta_2(2\la_1|2\tau)\theta_2(2\la_2|2\tau).
\end{align*}
Here $\la_i = \la_{\ep_i} \; (i=1,2)$.
Therefore we can prove this lemma 
by using the identities of theta functions 
(addition theorems) (\ref{44}), (\ref{33}), (\ref{22}), (\ref{11}) in
appendix.
\end{prf}

We note that the anti-symmetric function $\Theta_{\ep_1}(\la) - \Theta_{\ep_2}(\la)
= \theta_1(\la_+) \theta_1(\la_-)$ is also the eigenfunction with 
eigenvalue zero.

\section*{Acknowledgments}
The author is grateful to Gen Kuroki, Takeshi Ikeda,
Masatoshi Noumi, Yasuhiko Yamada, Yasushi Komori, and Kazuhiro Hikami
for fruitful discussions and kind interest.
He also thank Koichi Takemura who suggested the formula like (\ref{takemura})
based on the knowledge of the corresponding differential system,
and Koji Hasegawa for helpful conversations.

\section*{Appendix A: Differential limit}
\def\thesection{A}
\setcounter{equation}{0}

Let us clarify the connection between our system of difference operators
and a quantization of the Inozemtsev Hamiltonian \cite{IM} \cite{I2}. 
By expanding in $\h$ one infers that
\begin{align*}
\widetilde{M}_1 &= 4 + M_{1,2} \h^2 + M_{1,4} \h^4 + O(\h^5), \\
\widetilde{M}_2 &= 8 + M_{2,2} \h^2 + M_{2,4} \h^4 + O(\h^5).
\end{align*}
If we abbreviate a function $f(\la_{\ep_1 \pm \ep_2})$ as
$f(\pm)$, $\partial_i = \frac{\partial}{\partial \la_i}$ $(i=1,2)$,
and $\theta_1'(z) = \frac{d}{dz} \theta_1(z)$ etc. 
We have
\begin{align*}
M_{1,2}  = &
\partial_1^2+\partial_2^2 \\
&-2\left(
\frac{\theta'_1}{\theta_1}(+)
+\frac{\theta'_1}{\theta_1}(-)
\right)\partial_1 
-2\left(
\frac{\theta'_1}{\theta_1}(+)
-\frac{\theta'_1}{\theta_1}(-)
\right)\partial_2
\\
& +2\left(
\frac{\theta''_1}{\theta_1}(+)
+\frac{\theta''_1}{\theta_1}(-)
\right),\\
M_{2,2}  =&  2 M_{1,2},  
\end{align*}
and
\begin{eqnarray*}
M_{2,4} - 2 M_{1,4} &=& 
\partial_1^2\partial_2^2 \\
&& -2 \left(
\frac{\theta_1'}{\theta_1}(+)
-\frac{\theta_1'}{\theta_1}(-)
\right)\partial_1^2 \partial_2  
-2 \left(
\frac{\theta_1'}{\theta_1}(+)
+\frac{\theta_1'}{\theta_1}(-)
\right)\partial_1 \partial_2^2  \\
&& + \left\{ 2 \left(
\left(\frac{\theta_1'}{\theta_1}\right)^2(+)
+\left( \frac{\theta_1'}{\theta_1} \right)^2 (-)
\right) 
- \left(
\frac{\theta_1''}{\theta_1}(+)
+2\frac{\theta_1'}{\theta_1}(+)\frac{\theta_1'}{\theta_1}(-)
+\frac{\theta_1''}{\theta_1}(-)
\right) \right\} \partial_1^2 \\
&& + \left\{ 2 \left(
\left(\frac{\theta_1'}{\theta_1}\right)^2(+)
+\left( \frac{\theta_1'}{\theta_1} \right)^2 (-)
\right) 
- \left(
\frac{\theta_1''}{\theta_1}(+)
-2\frac{\theta_1'}{\theta_1}(+)\frac{\theta_1'}{\theta_1}(-)
+\frac{\theta_1''}{\theta_1}(-)
\right) \right\} \partial_2^2 \\
&& +  4 \left(
\left(\frac{\theta_1'}{\theta_1}\right)^2(+)
- \left( \frac{\theta_1'}{\theta_1} \right)^2 (-)
\right)  \partial_1 \partial_2 \\
&& + \left\{ 2 \left(
\frac{\theta_1'\theta_1''}{\theta_1^2}(+)
+\frac{\theta_1'\theta_1''}{\theta_1^2}(-)
\right)  
+2 \left(
\frac{\theta_1''}{\theta_1}(+)\frac{\theta_1'}{\theta_1}(-)
+\frac{\theta_1'}{\theta_1}(+)\frac{\theta_1''}{\theta_1}(-)
\right)
\right. \\
&&\left. {}
-4 \left( \left( \frac{\theta_1'}{\theta_1} \right)^3 (+)
+\left( \frac{\theta_1'}{\theta_1} \right)^3(-)
\right) \right\} \partial_1 \\
&& + \left\{ 2 \left(
\frac{\theta_1'\theta_1''}{\theta_1^2}(+)
-\frac{\theta_1'\theta_1''}{\theta_1^2}(-)
\right)  
-2 \left(
\frac{\theta_1''}{\theta_1}(+)\frac{\theta_1'}{\theta_1}(-)
-\frac{\theta_1'}{\theta_1}(+)\frac{\theta_1''}{\theta_1}(-)
\right)
\right. \\
&&\left. {} -4 \left( \left( \frac{\theta_1'}{\theta_1} \right)^3 (+)
-\left( \frac{\theta_1'}{\theta_1} \right)^3(-)
\right) \right\} \partial_2 \\
&&+ \frac{1}{2}
\left( \frac{\theta_1^{(4)}}{\theta_1}(+) + \frac{\theta_1^{(4)}}{\theta_1}(-) \right)
 -4 
\left( \frac{\theta_1'''\theta_1'}{\theta_1^2}(+) + \frac{\theta_1'''\theta_1'}{\theta_1^2}(-) \right) \\
&& + 2 \left(\frac{\theta_1''\theta_1'{}^2}{\theta_1^3}(+) +
\frac{\theta_1''\theta_1'{}^2}{\theta_1^3}(-) \right) 
-2 \frac{\theta_1''}{\theta_1}(+)\frac{\theta_1''}{\theta_1}(-),  
\end{eqnarray*}

We set $\Delta = \theta_1(+) \theta_1(-)$, then 
\begin{eqnarray}
\Delta^{-1} \cdot M_{2,2} \cdot \Delta
&=& \partial_1^2 + \partial_2^2 
+ 4 \left( \left(\frac{\theta_1''}{\theta_1} - \frac{\theta_1'{}^2}{\theta_1^2}\right)(+)
 + \left( \frac{\theta_1''}{\theta_1} - \frac{\theta_1'{}^2}{\theta_1^2} \right)(-) \right) \nonumber \\
&=& \partial_1^2 + \partial_2^2 
+ 4 \left( (\log\theta_1)''(+) + (\log\theta_1)''(-) \right),  \label{bibun}
\end{eqnarray}
\begin{eqnarray*} 
&& \Delta^{-1} \cdot (M_{2,4} - 2M_{1,4}) \cdot \Delta \\
&=& \partial_1^2\partial_2^2  \\
&&+ 4 \left( \left( \frac{\theta_1''}{\theta_1}+ \frac{\theta_1'{}^2}{\theta_1^2}\right) (+) 
- \left( \frac{\theta_1''}{\theta_1}+ \frac{\theta_1'{}^2}{\theta_1^2}\right) (-) \right) 
\partial_1\partial_2 \\
&& +2 \left(  \left( \frac{\theta_1'''}{\theta_1} - 3 \frac{\theta_1'' \theta_1'}{\theta_1^2}
+2\frac{\theta_1'{}^3}{\theta_1^3}\right) (+)
+ \left( \frac{\theta_1'''}{\theta_1} - 3 \frac{\theta_1'' \theta_1'}{\theta_1^2}
+2\frac{\theta_1'{}^3}{\theta_1^3}\right) (-) \right) \partial_1 \\
&& + 2 \left(  \left( \frac{\theta_1'''}{\theta_1} - 3 \frac{\theta_1'' \theta_1'}{\theta_1^2}
+2\frac{\theta_1'{}^3}{\theta_1^3}\right) (+)
- \left( \frac{\theta_1'''}{\theta_1} - 3 \frac{\theta_1'' \theta_1'}{\theta_1^2}
+2\frac{\theta_1'{}^3}{\theta_1^3}\right) (-) \right) \partial_2 \\
&&+ 2 \left( \frac{\theta_1^{(4)}}{\theta_1}(+) + \frac{\theta_1^{(4)}}{\theta_1}(-) \right)
-8 \left( \frac{\theta_1''' \theta_1'}{\theta_1^2}(+) + \frac{\theta_1''' \theta_1'}{\theta_1^2}(-) \right) \\
&& -2 \left( \left( \frac{\theta_1''}{\theta_1} \right)^2 (+) +
 \left( \frac{\theta_1''}{\theta_1} \right)^2 (-) \right)
-8  \frac{\theta_1''}{\theta_1}(+)  \frac{\theta_1''}{\theta_1}(-) \\
&& +16 \left(\frac{\theta_1'' \theta_1'{}^2}{\theta_1^3}(+) + \frac{\theta_1'' \theta_1'{}^2}{\theta_1^3}(-) \right) 
+8 \left( \frac{\theta_1'{}^2}{\theta_1^2}(+)\frac{\theta_1''}{\theta_1}(-) + 
\frac{\theta_1''}{\theta_1}(+)\frac{\theta_1'{}^2}{\theta_1^2}(-) \right) \\
&&-8 \left( \left( \frac{\theta_1'}{\theta_1} \right)^4 (+) +
 \left( \frac{\theta_1'}{\theta_1} \right)^4 (-)  \right)
 +8 \left(\frac{\theta_1'}{\theta_1} \right)^2(+) \left( \frac{\theta_1'}{\theta_1}\right)^2(-) \\
&=& \partial_1^2\partial_2^2  \\
&&+4 \{ (\log\theta_1)''(+) -
(\log\theta_1)''(-) \} 
\partial_1\partial_2  
\\
&& + 2 \left\{ (\log\theta_1)'''(+) +
(\log\theta_1)'''(-) \right\} 
\partial_1
+ 2 \{ (\log\theta_1)'''(+) -
(\log\theta_1)'''(-) \} 
\partial_2
\\
&& + 2 \{ (\log\theta_1)^{(4)}(+) + (\log\theta_1)^{(4)}(-) \} \\
&& + 4 \{ (\log\theta_1)''(+) - (\log\theta_1)''(-) \}^2 \\
&=& \left\{ \partial_1 \partial_2 + 2 \left( (\log \theta_1)''(+) -
(\log \theta_1)''(-) \right) \right\}^2
\end{eqnarray*}

The complete integrable Hamiltonian of type $BC_n$ is introduced by
Olshanetsky-Perelomov \cite{OP}, and later generated by
Inozemtsev-Meshcheryakov \cite{IM} \cite{I2}.
In the rank two case, the Hamiltonian is
\begin{align*}
H = &-\frac{1}{2}(\partial_1^2+\partial_2^2) 
+ g(g-1) \left( \wp(x_1+ x_2) + \wp(x_1- x_2) \right)  \\
&+ \sum_{0 \leq r \leq 3}  g_r (g_r-1)
 (\wp(\omega_r + x_1) +  \wp(\omega_r + x_2)),
\end{align*}
where $\wp (x)$ denotes the Weierstrass $\wp$-function
with two periods $2\omega_1$ and $2\omega_2$,
and $\omega_0 = 0, \; \omega_3 = -\omega_1 - \omega_2$.
By the connection between theta function and $\wp$ function (\ref{pekan})
in Appendix B, our differential limit (\ref{bibun}) is identified with this
Hamiltonian for the special coupling 
constants $g (g-1) =2$, and $g_r(g_r-1) = 0\; (0 \leq r \leq 3)$.

\section*{Appendix B: Theta function}
\def\thesection{B}
\setcounter{equation}{0}

We establish notations and identities on the theta functions \cite{WW}.
The Jacobi theta functions are defined for $\tau \in \mathfrak{H}_+$ as
follows:
\[ \theta_1(z|\tau) = \sum_{k \in \itg}
\exp\left[ 2\pi i \left( (z+\frac{1}{2})(k+\frac{1}{2}) +
\frac{1}{2} (k+\frac{1}{2})^2 \tau \right) \right] \] 
\[ \theta_2(z|\tau) = \sum_{k \in \itg}
\exp\left[ 2\pi i \left(  z(k+\frac{1}{2})  + \frac{1}{2} (k+\frac{1}{2})^2 \tau 
\right) \right] \] 
\[ \theta_3(z|\tau) = \sum_{k \in \itg}
\exp \left[ 2 \pi i \left( zk + \frac{k^2}{2}
\tau \right) \right] \] 
\[ \theta_4(z|\tau) = \sum_{k \in \itg}
\exp\left[ 2\pi i \left( (z+\frac{1}{2})k + \frac{k^2}{2} \tau \right) 
\right] \] 
Note that $\theta_1(z)$ is odd and the
other three are even.
These functions has quasi-periodicity:
\begin{align}
\theta_1(z+m|\tau) &=(-1)^m \theta_1(z|\tau), \label{henkan1} \\
\theta_1(z+m\tau|\tau) &=(-1)^me^{-\pi im^2\tau-2\pi imz}\theta_1(z|\tau),
\label{henkan2} 
\end{align}
$(m\in \itg)$,
while other three can be expressed by $\theta_1(z)$
\begin{align} 
\theta_1(z+\frac{1}{2}|\tau) &=  \theta_2(z|\tau), \nonumber \\ 
\theta_1(z+\frac{\tau}{2}|\tau) 
&= i e^{-\pi i \left(z+\frac{\tau}{4} \right)} \theta_4(z|\tau), \label{theta2} \\
\theta_1(z+\frac{1}{2}+\frac{\tau}{2} |\tau) 
&= e^{-\pi i \left(z+\frac{\tau}{4} \right)} \theta_3(z|\tau). \nonumber
\end{align}

We use these identities in the computations in Lemma \ref{kitei}.
\begin{align}
\theta_4(x|\tau) \theta_4(y|\tau) &= 
\theta_3(x+y|2\tau) \theta_3(x-y|2\tau) -
\theta_2(x+y|2\tau) \theta_2(x-y|2\tau), \label{44} \\
\theta_3(x|\tau) \theta_3(y|\tau) &= 
\theta_3(x+y|2\tau) \theta_3(x-y|2\tau) +
\theta_2(x+y|2\tau) \theta_2(x-y|2\tau), \label{33} \\
\theta_2(x|\tau) \theta_2(y|\tau) &= 
\theta_3(x+y|2\tau) \theta_2(x-y|2\tau) +
\theta_2(x+y|2\tau) \theta_3(x-y|2\tau), \label{22} \\
\theta_1(x|\tau) \theta_1(y|\tau) &= 
\theta_3(x+y|2\tau) \theta_2(x-y|2\tau) -
\theta_2(x+y|2\tau) \theta_3(x-y|2\tau). \label{11}
\end{align}

The sigma function $\sigma(z)$ is an entire, odd, and quasi-periodic
function with two primitive quasi-periods $2\omega_1, 2\omega_2$.
\[ \sigma(z + 2n\omega_1 + 2m\omega_2) 
= (-1)^{n+m+nm} e^{(2n\eta_1 + 2m\eta_2)(z + n\omega_1 + m\omega_2)} \sigma(z) \]
with $\eta_i = \zeta(\omega_i)$ $(i=1,2)$, where $\zeta(z) = \sigma'(z)/\sigma(z)$ 
denotes the Weierstrass $\zeta$-function.
The connection between the Jacobi theta functions and the
sigma functions are
\[ \sigma(z) = \left( \exp \frac{\eta_1 z^2}{2 \omega_1} \right)
\frac{\theta_1(z/2\omega_1)}{\theta_1'(0)},  \]
\[
\quad \sigma_r(z) = \left( \exp \frac{\eta_1 z^2}{2 \omega_1} \right)
\frac{\theta_{r+1}(z/2\omega_1)}{\theta_{r+1}(0)} \;
(r = 1,2,3).\]
Then, for the function $v(z)$ in van Diejen's system (\ref{vz}), we have
\begin{equation}\label{tesig}
v(z) := \frac{\sigma(z+\mu)}{\sigma(z)} =
\left( \exp \frac{\eta_1 (2z\mu + \mu^2)}{2 \omega_1} \right)
\frac{\theta_1((z+\mu)/2\omega_1)}{\theta_1(z/2\omega_1)}.
\end{equation}
The connection with $\wp$ function is
\begin{equation}\label{pekan}
 \wp(z) = - \frac{d^2}{dz^2} \log \sigma(z) = 
 - \frac{1}{4\omega_1^2} \left( \frac{d^2}{dz^2} \log \theta_1(z/2\omega_1) \right)
- \frac{\eta_1}{\omega_1}.
\end{equation}


\begin{thebibliography}{BKMS}

\bibitem
[Bax]{Baxter}
R.J. Baxter,
``Eight-vertex model in lattice statistics and
one-dimensional anisotropic Heisenberg chain.''
I. Ann.\ Phys. {\bf 76}(1973), 1-24,
II. {\it ibid.} 25-47,
III. {\it ibid.} 48-71.


\bibitem
[Ch]{C}
Cherednik, Ivan 
``Double affine Hecke algebras and Macdonald's conjectures'',
Ann. of Math. (2) {\bf 141}(1995), 191-216.


\bibitem
[vD1]{vD2}
J.F.van Diejen,
``Integrability of difference Calogero-Moser systems",
J. Math. Phys. {\bf 35}(1994), 2983-3004.

\bibitem
[vD2]{vD1}
J.F. van Diejen,
``Commuting difference operators with polynomial eigenfunctions",
Compositio Math. {\bf 95}(1995), 183-233.


\bibitem
[FV1]{FV1}
G.Felder, A.Varchenko, 
``Algebraic Bethe ansatz for the elliptic quantum group $E_{\tau,\eta}(sl_2)$,''  
Nucl. Phys. B {\bf 480}(1996) 485--503.


\bibitem
[FV2]{FV2}
G.Felder, A.Varchenko, ``Elliptic
quantum groups and Ruijsenaars models,'' 
J. Statist. Phys. {\bf 89}(1997),
963-980.


\bibitem
[H1] {has}
K.Hasegawa,
``On the crossing symmetry of the
elliptic solution of the Yang-Baxter equation
and a new L operator for Belavin's solution'',
J.Phys. A: Math. Gen. {\bf 26} (1993), 3211-3228.

\bibitem
[H2]{has93}
K. Hasegawa,
``L-operator for Belavin's R-matrix acting on the space of theta
functions'',
J.Math.Phys. {\bf 35}(1994),  6158-6171.


\bibitem [H3]      
{has97}
K.Hasegawa
``Ruijsenaars' Commuting difference operators as commuting transfer 
matrices'',
Commun. Math. Phys. {\bf 187} (1997), 289-325.


\bibitem[HIK]{HIK}
K.Hasegawa, T.Ikeda, T.Kikuchi, 
``Commuting difference operators arising from 
the elliptic $C_2^{(1)}$-face model'',
J. Math. Phys. {\bf 40} (1999), 4549-4568.


\bibitem[IM]{IM}
V.I.Inozemtsev, D.V.Meshcheryakov, 
``Extension of the class of integrable dynamical systems connected with 
semisimple Lie algebras''
Lett.\ Math.\ Phys.
{\bf 9} (1985), 13-18.


\bibitem[I]{I2}
V.I.Inozemtsev,
``Lax representation with spectral parameter on a torus for integrable
particle systems",
Lett.\ Math.\ Phys.
{\bf 17} (1989), 11-17.


\bibitem[JMO1]{JMO88}
M. Jimbo, T. Miwa and M. Okado,
``Local
state probabilities of solvable lattice models:
An $A_n^{(1)}$ family'',
Nucl.\ Phys. B{\bf 300}(1988), 74-108.


\bibitem[JMO2]{JMO2}
M. Jimbo, T. Miwa and M. Okado,
``Solvable lattice models related to the vector representation
of classical simple Lie algebras",
Commun.\ Math.\ Phys.
{\bf 116} (1988), 507-525.

\bibitem
[KP]{KP}
V.G.Kac and D.H.Peterson
``Infinite-dimensional Lie algebras, theta functions and modular forms'',
Adv. in Math. {\bf 53} (1984), 125-264


\bibitem
[KH1]{KH97}
Y. Komori, K. Hikami,
``Quantum integrability of the generalized elliptic Ruijsenaars
models'',
J. Phys. A: Math. Gen. {\bf 30}(1977), 4341-4364.

\bibitem
[KH2]{KH98}
Y. Komori, K. Hikami,
``Conserved operators of the generalized elliptic Ruijsenaars
models'',
J. Math. Phys. {\bf 39}(1998) 6175-6190 .

\bibitem
[KH3]{KH_K-op1}
Y. Komori and K. Hikami,
``Notes on operator-valued solutions of the Yang-Baxter equation and
the reflection equation'',
Mod. Phys. Lett. A {\bf 11}(1996), 2861-2870.

\bibitem
[KH4]{KH_K-op2}
Y. Komori and K. Hikami,
``Elliptic $K$-matrix associated with Belavin's symmetric
$R$-matrix'',
Nucl. Phys. B {\bf 494}(1997), 687-701.



\bibitem
[K]
{Koorn}
T.H.Koornwinder,
``Askey-Wilson polynomials for root systems of type $BC$'',
Contemp.Math. {\bf 138} (1992), 189-204. 

\bibitem
[M]{Macbook}
I.G.Macdonald,
{\it Symmetric functions and Hall polynomials}(2nd ed.),
Oxford Univ. Press, 1995.

\bibitem[N]{N}
M. Noumi, 
``Macdonald-Koornwinder polynomials and affine Hecke rings
(Japanese)'',
Various aspects of hypergeometric functions (Japanese) (Kyoto, 1994).
S\=urikaisekikenky\=usho K\=oky\=uroku No. 919(1995),
44-55.


\bibitem
[OP]{OP}
M.A.Olshanetsky and A.M.Perelomov,
``Completely integrable Hamiltonian systems connected with semi-simple
Lie algebras'', Inv. Math. {\bf 37}(1976), 93-108.


\bibitem
[R]{Ruij}
S.N.M.Ruijsenaars,
``Complete integrability of relativistic Calogero-Moser systems
and elliptic function identities", Comm.Math.Phys. {\bf 110} (1987), 191-213.


\bibitem
[SU]{SU}
Y. Shibukawa and K. Ueno, 
``Completely ${\mathbb Z}$ symmetric $R$-matrix'',
Lett. Math. Phys. {\bf 25}(1992), 239-248.


\bibitem
[WW]{WW}
E.T.Whittaker, G.N. Watson,
``A course of modern analysis'',
Cambridge: Cambridge U. P. 1986

\end{thebibliography}
\end{document}